\documentclass{amsart}
\usepackage{amssymb}

\usepackage[alphabetic]{amsrefs}
\usepackage[all]{xy}
\usepackage{tikz}
\usetikzlibrary{matrix}

\newtheorem{lemma}{Lemma}[section]

\newtheorem{proposition}[lemma]{Proposition}
\newtheorem{theorem}[lemma]{Theorem}

\newtheorem{corollary}[lemma]{Corollary}

\theoremstyle{definition}   
\newtheorem{remark}[lemma]{Remark}
\newtheorem{definition}[lemma]{Definition}
\newtheorem{example}[lemma]{Example}

\def\<{\langle}
\def\>{\rangle}
\def\0{\<0\>}
\def\BL{\mathsf{BL}}
\def\BA{\mathsf{BA}}
\def\DL{\mathsf{DL}}
\def\Spec{{\sf Spec\;}}
\def\Mod{\mbox{\sf Mod}}
\def\Hom{\mbox{Hom}}
\def\bz{\mathbb{Z}}
\def\bq{{\mathbb{Q}}}
\def\Fp{\mathbb{F}_p}
\def\zp{\bz_{(p)}}
\def\cC{\textsf{T}}
\def\cB{\textsf{S}}
\def\S{\mathcal{S}}
\def\K{\textsf{K}}
\def\Lat{\mathsf{L}}
\def\loc{{\sf loc}}
\def\th{{\sf th}}
\def\max{\textsf{Max}}
\def\w{\wedge}
\def\da{{\downarrow}}
\def\ua{{\uparrow}}
\def\lra{\longrightarrow}
\def\modJ{\;(\mbox{mod}\; J)}

\def\Xc{X^c}
\def\Ll{\Lambda}
\def\G{{\Lambda_{\zp}}}
\def\Pi{{\Lambda_{\bq}}}
\def\L{{\Lambda_{\Fp}}}
\def\IN{I\Ll}
\def\f{f_\bullet}
\def\i{f^\bullet}
\def\g{g_\bullet}
\def\ig{g^\bullet}
\def\h{h_\bullet}
\def\ih{h^\bullet}
\def\oG{\overline{\f}}
\def\DRiS{ \BL\left( D(R)/\<\i S\>\right)}

\newcommand{\one}{{\rm 1\hspace*{-0.4ex}
\rule{0.1ex}{1.52ex}\hspace*{0.2ex}}}

\title[Bousfield lattices of non-Noetherian rings]{Bousfield lattices of non-Noetherian rings: \newline some quotients and products}

\begin{document}

%
%

\author{F. Luke Wolcott}             

\email{luke.wolcott@gmail.com}

\address{Department of Mathematics,
         Lawrence University,
         Appleton, WI, 54911,
         USA}



\begin{abstract}
In the context of a well generated tensor triangulated category, Section 3 investigates the relationship between the Bousfield lattice of a quotient and quotients of the Bousfield lattice.  In Section 4 we develop a general framework to study the Bousfield lattice of the derived category of a commutative or graded-commutative ring, using derived functors induced by extension of scalars.  Section 5 applies this work to extend results of Dwyer and Palmieri to new non-Noetherian rings.
\end{abstract}


\maketitle

\section{Introduction} \label{section-intro}

Let $R$ be a commutative ring and consider  the unbounded derived category $D(R)$ of right $R$-modules.  Given an object $X\in D(R)$, define the Bousfield class $\<X\>$ of $X$ to be $\{W\in D(R)\;|\; W\otimes_R^L X=0\}$.  Order Bousfield classes by reverse inclusion, so $\<0\>$ is the minimum and $\<R\>$ is the maximum.  It is known that there is a set of such Bousfield classes.  The join of any set $\{\<X_\alpha\>\}$ is the class $\<\coprod_\alpha X_\alpha\>$, and the meet of a set of classes is the join of all the lower bounds.  The collection of Bousfield classes thus forms a lattice, called the Bousfield lattice $\BL(D(R))$.

A full subcategory of $D(R)$ is localizing if it is closed under triangles and arbitrary coproducts.  Thus every Bousfield class is a localizing subcategory.  A result of Neeman's~\cite{[Nee]} shows that when $R$ is Noetherian, every localizing subcategory is a Bousfield lattice, and this lattice is isomorphic to the lattice of subsets of the prime spectrum $\Spec R$.

The case of a non-Noetherian ring is much less understood.  Given a ring $k$, fix $n_i>1$ and define
\[\Lambda_k = \frac{k[x_1, x_2, x_3,...]}{(x_1^{n_1}, x_2^{n_2}, x_3^{n_3},...)}, \] and give $\Lambda_k$ a grading by setting $deg(x_i)=2^i$.  Consider the unbounded derived category $D(\Lambda_k)$ of right graded $\Lambda_k$-modules; objects in $D(\Ll_k)$ are bi-graded.  Dwyer and Palmieri~\cite{[DP]} studied the Bousfield lattice of this category, when $k$ is a countable field (see Example~\ref{example-lambda} below for more details).  The initial motivation for the present work was to extend their main results to the case where $k=\zp$.  We have done this fairly completely, and in the process developed tools that apply in much more general settings.

Our most general results apply to the Bousfield lattice of any well generated tensor triangulated category, and appear in Section~\ref{section-landq}.  Iyengar and Krause~\cite{[IK]} recently showed that a well generated tensor triangulated category has a set of Bousfield classes, and thus a Bousfield lattice.  Note that compactly generated tensor triangulated categories are well generated, and in particular those generated by the tensor unit are.  For simplicity in this introduction, suppose $\cC$ is a tensor triangulated category generated by the tensor unit $\one$; this includes the case of the derived category of a ring, but also the stable homotopy category and the stable module category of a $p$-group.  Let $-\w-$ denote the tensor product, and $-\vee-$ denote the join.

The results of Section~\ref{section-landq} concern the relationship between the quotient of a lattice and the lattice of a quotient.  Given $Z\in \cC$, consider the Verdier quotient $\cC/\<Z\>$; this is well generated because $\<Z\>$ is.  The quotient functor $\pi\colon\cC\to \cC/\<Z\>$ induces a well-defined, order-preserving map of lattices $\pi\colon\BL(\cC)\to \BL(\cC/\<Z\>)$, where $\<X\>\mapsto \<\pi X\>$.

Given $Z\in \cC$, define $a\<Z\>$ to be the join of all classes $\<Y\>$ such that $\<Z\w Y\>=\0$.  For any class $\<X\> \in \BL(\cC)$, define $\<X\>\da$ to be the collection of classes less than or equal to $\<X\>$.  In Definition~\ref{defn-quotient} we give a notion of quotient lattice.  The following is Proposition~\ref{epim}.  


\begin{proposition} Let $\<Z\>$ be any Bousfield class in $\BL(\cC)$.  Then $\pi$ induces an onto lattice join-morphism with trivial kernel
\[ \overline{\pi}\colon \BL(\cC)/(a\<Z\>)\da\lra \BL\left(\cC/\<Z\>\right).\] 
\end{proposition}

In this context, we say a class $\<X\>$ is complemented if $\<X\>\vee a\<X\>=\<\one\>$.  The sub-poset of complemented classes is denoted $\BA(\cC)$.


\begin{corollary} If $\<Z\>$ is complemented then the above map is an isomorphism of lattices.
\end{corollary}

This is proven in Corollary~\ref{compl}.  We also consider the sub-poset $\DL(\cC)$ of classes $\<X\>$ such that $\<X\w X\> = \<X\>$.  The following is Proposition~\ref{not-isomo}.


\begin{proposition} If $\<Z\>$ is an element of $\DL(\cC)$ but is not complemented, then the map in Proposition~\ref{epim} is not an isomorphism.  This happens in the stable homotopy category and in $D(\Lambda_k)$, where $k$ is a countable field.
\end{proposition}

These results rely in part on an interesting observation that we have been unable to find in the literature.  Call an object $X\in \cC$ square-zero if $X$ is nonzero but $X\w X=0$.  In Corollary~\ref{sq-free-3} we prove the following.


\begin{corollary}There are no square-zero objects in $\cC$ if and only if
\[\BL = \DL = \BA.\]
\end{corollary}

Section~\ref{section-rmbl} specializes to look at functors between derived categories of rings.  A ring map $f\colon R\to S$ induces a functor $f_*\colon \Mod$-$R \to \Mod$-$S$, via extension of scalars, and the forgetful functor $f^*$ is a right adjoint.  This carries to the level of chain complexes, and we get an adjoint pair of derived functors on derived categories \[\f\colon D(R)\rightleftarrows D(S)\colon \i.\] 

These functors induce maps between lattices, where $\f\<X\>=\<\f X\>$ and $\i\<Y\> = \<\i Y\>$, which preserve order and arbitrary joins.  

First we investigate the behavior of the sub-posets $\BA$ and $\DL$ under $\f$ and $\i$.  Let $\<M_f\>$ be the join of all classes $\<Y\>$ with $\f\<Y\>=\0$.  Abbreviate $\BL(D(R))$ to $\BL(R)$ or $\BL_R$, and likewise for $\BA$ and $\DL$.  Our most general statement is the following, which is Proposition~\ref{G_i_sl}.  


\begin{proposition} Suppose $\f \i\<X\>=\<X\>$ for all $\<X\>$.  The following hold.
\begin{enumerate}
    \item \textit{The map $\f $ sends $\DL_{R}$ onto $\DL_{S}$, and the map $\i$ injects $\DL_{S}$ into $\DL_{R}$.}
    \item \textit{The map $\f$ sends $\BA_{R}$ onto $\BA_{S}$, and if $\<\i S\>\vee \<M_f\>=\<R\>$ then  $\i$ injects $\BA_{S}$ into $\BA_{R}$.}
\end{enumerate}    
\end{proposition}

Next we establish maps between various quotients and lattices.  

\begin{center}
\begin{tikzpicture}
  \matrix (m) [matrix of math nodes, row sep=3em,
    column sep=3em,ampersand replacement=\&]{
       \BL(R) \& \BL(S) \\
       \BL(R)\slash \<M_f\>\da \& \BL(D(R)/\<\i S\>)  \\}; 
    \path[>=latex,->,font=\scriptsize]
       (m-1-1) edge (m-2-1)
       (m-1-1) edge node[auto]{$\f$} (m-1-2)
       (m-2-1) edge node[auto,swap]{$(*)$} (m-2-2)
       (m-2-2) edge node[auto,swap]{$(\dagger)$} (m-1-2)
       (m-2-1) edge node[auto]{$\exists$} (m-1-2);
   
\end{tikzpicture}
\end{center}

We show $\<M_f\> = a\<\i S\>$, and so Corollary~\ref{compl} implies that the map $(*)$ is an isomorphism when $\<\i S\>\vee\<M_f\>=\<R\>$.  Theorem~\ref{latt-o-q} states that the map $(\dagger)$ exists and is an isomorphism when $\f \i\<X\>=\<X\>$ for all $\<X\>$.

Finally, Section~\ref{section-nnr} applies the results of the previous two sections.  Let $g\colon  \Ll_{\zp}\to \Ll_{\Fp}$ be the obvious projection, and let $h\colon  \Ll_{\zp}\to \Ll_{\bq}$ be inclusion.  These maps give derived functors and lattice maps, as above.  The map $\g$ has $\<\g\ig X\>=\<X\>$ for all $\<X\>$ (Proposition~\ref{fiX_X_gamma}), and so the diagram above becomes the following.

\begin{center}
\begin{tikzpicture}
  \matrix (m) [matrix of math nodes, row sep=3em,
    column sep=3em,ampersand replacement=\&]{
       \BL(\G) \& \BL(\L) \\
       \BL(\G)\slash \<\ih\Pi\>\da \& \BL(D(\G)/\<\ig \L\>)  \\}; 
    \path[>=latex,->,font=\scriptsize]
       (m-1-1) edge (m-2-1)
       (m-1-1) edge node[auto]{$\g$} (m-1-2)
       (m-2-1) edge node[auto,swap]{$\cong$} (m-2-2)
       (m-2-2) edge node[auto,swap]{$\cong$} (m-1-2)
       (m-2-1) edge (m-1-2);
   
\end{tikzpicture}
\end{center}

Theorem~\ref{split-1} gives a splitting of the Bousfield lattice of $\G$ as the product lattice $\<\ig\L\>\da\times \<\ih\Pi\>\da$.  Combining this with other results, we conclude the following. Let $\loc(X)$ denote the smallest localizing subcategory containing $X$.


\begin{corollary} The functors $\g$ and $\h$ induce lattice isomorphisms
\[ \BL({\G}) \cong \BL({\L})\times \BL({\loc(\ih\Pi)}),\]  
\[ \DL({\G}) \cong \DL({\L})\times \DL({\loc(\ih\Pi)}),\]
\[ \BA({\G}) \cong \BA({\L})\times \BA({\loc(\ih\Pi)}), \]
\[ \mbox{ where } \<X\> \mapsto \left(\g\<X\>, \<X\w\ih\Pi\>\right).\]
\end{corollary}

This is proven in Corollaries~\ref{split-2} and \ref{split-3}.  As immediate corollaries to this, we get that the cardinality of $\BL(\G)$ is $2^{2^{\aleph_0}}$ (Corollary~\ref{size}) and that, unlike in $\BL(\L)$, in $\BL(\G)$ there is no nonzero minimum Bousfield class (Proposition~\ref{no-min}).

Section~\ref{section-BLbackground} contains background on Bousfield lattices and gives examples.  With the exception of Proposition~\ref{sq-free-1} and its corollaries, and our treatment of Bousfield lattices of proper subcategories, the contents are not new.  The results of the Sections$~\ref{section-landq} - \ref{section-nnr}$ are new, unless cited.  

We are grateful to John Palmieri and Dan Christensen for many helpful conversations and suggestions.

\section{Background on Bousfield lattices} \label{section-BLbackground}

In this section we review the definition and basic properties of Bousfield classes and the Bousfield lattice, and outline some of what is known about the Bousfield lattice in several examples.  Most of the following general properties of Bousfield classes were first established by Bousfield~\cites{[Bou79a],[Bou79b]} in the context of the stable homotopy category.  Further work was done in~\cites{[ravenel],[HPS],[HP], [IK]}.  Our lattice theory reference is~\cite{[Bir]}.  We will work in the context of a well generated tensor triangulated category, which we now define.\\

\begin{definition}~\cite[\textsection 6.3]{[krause-loc-survey]} Let $\cC$ be a triangulated category which admits arbitrary coproducts and fix a regular cardinal $\alpha$.  An object $X$ in $\cC$ is called \textit{$\alpha$-small} if every morphism $X\to \coprod_{i\in I} Y_i$ in $\cC$ factors through $\coprod_{i\in J}Y_i$ for some subset $J\subseteq I$ with $card(J) < \alpha$.  The triangulated category $\cC$ is called \textit{$\alpha$-well generated} if it is perfectly generated by a set of $\alpha$-small objects (see~\cite[\textsection 5.1]{[krause-loc-survey]}).  And $\cC$ is called \textit{well generated} if it is $\beta$-well generated for some regulard cardinal $\beta$.
\end{definition}

A category is $\aleph_0$-well generated if and only if it is compactly generated.  A triangulated category is \textit{tensor triangulated} if it has a symmetric monoidal product, which we will denote $-\w-$, that is compatible with the triangulation, is exact in both variables, and commutes with arbitrary coproducts~\cite[App.~A]{[HPS]}.  Let $\Sigma$ denote the shift.  We will denote the tensor unit by $\one$, and do not assume that $\one$ is compact.  

\begin{definition} Let $X$ be an object in $\cC$.
\begin{enumerate}
   \item A full subcategory $\cB\subseteq \cC$ is \textit{thick} if it is closed under triangles and retracts.
   \item The smallest thick subcategory containing $X$ is denoted $\th(X)$; this is also called the \textit{thick subcategory generated by $X$}.
   \item A full subcategory $\cB\subseteq \cC$ is \textit{localizing} if it is closed under triangles, retracts, and arbitrary coproducts.
   \item A full subcategory $\cB\subseteq \cC$ is a \textit{tensor ideal} if $X\in \cB$ and $Y\in \cC$ implies $X\w Y\in \cB$.
   \item The smallest localizing subcategory containing $X$ is denoted $\loc(X)$; this is also called the \textit{localizing subcategory generated by $X$.}
\end{enumerate}
\end{definition}

Note that if $\cC=\loc(\one)$, then every localizing subcategory $\cB\subseteq \cC$ is a tensor ideal.  Indeed, for $X\in \cB$ and $Y\in \loc(\one)$, then $X\w Y\in \loc(X\w\one)=\loc(X)\subseteq \cB$.

Henceforth, let $\cC$ denote a well generated tensor triangulated category, or a well generated localizing tensor ideal of a tensor triangulated category.  In the former case, of course we have $\one\in \cC$.  However, in the latter case we may have $\one\notin\cC$, and this introduces new subtleties in the structure of the Bousfield lattice.  

\begin{definition} Let $W, X$, and $Y$ be objects in $\cC$.
\begin{enumerate}
  \item  We say $W$ is \textit{$X$-acyclic} if $W\w X=0$.
  \item The collection of $X$-acyclics is denoted $\<X\>$ and called the \textit{Bousfield class} of $X$.
  \item We say $X$ and $Y$ are \textit{Bousfield equivalent} if they have the same acyclics.
\end{enumerate}
\end{definition}

There is a partial ordering on Bousfield classes, given by reverse inclusion.  So we say
\[\<X\>\leq \<Y\> \mbox{ if and only if } \left(W\w Y=0 \implies W\w X =0\right).\]

Note that $\0$ is the minimum class under this ordering.  When $\one\in\cC$, then $\<X\>=\0$ implies $X=0$.  The join of a set of classes $\{\<X_\alpha\>\}_{\alpha\in A}$ is given by 
\[\bigvee_{\alpha\in A} \<X_\alpha\> = \left\< \coprod_{\alpha\in A} X_\alpha \right\>.\]

It was recently shown in~\cite[Thm.~3.1]{[IK]} that in a well generated tensor triangulated category there is always a set of Bousfield classes.  Their proof applies as well to the setting of a well generated localizing tensor ideal of a tensor triangulated category.  We can define the meet (denoted $\curlywedge$) of any set of classes $\{\<X_\alpha\>\}$ to be the join of all the lower bounds; this join is over a set, and a nonempty set because $\0$ is the minimum.

A partially ordered set with finite joins and meets is called a \textit{lattice}.  A lattice with arbitrary joins and meets is \textit{complete}.  The collection of Bousfield classes of $\cC$ is thus a complete lattice, called the \textit{Bousfield lattice}, and denoted $\BL$.  

In any complete lattice there is also a maximum element $\<\max\>$, given by joining all elements.  When $\one\in \cC$, then clearly $\<\max\> = \<\one\> = \{0\}$.  On the other hand, see Remark~\ref{max-x}.

Given any well generated localizing tensor ideal $\cB \subseteq \cC$, we can consider the Bousfield lattice $\BL(\cB)$.  Some care is necessary, since for $X\in \cB$, the Bousfield class $\<X\>$ in $\BL(\cB)$ is $\{ W\in \cB\;|\; X\w W=0\}$.  If $X,Y\in \cB$ have $\<X\>\leq \<Y\>$ in $\BL(\cB)$, it does not necessarily follow that $\<X\>\leq \<Y\>$ in $\BL(\cC)$.  However, see Lemmas~\ref{incl-gh} and \ref{incl-sub}.

The tensor product gives another operation on Bousfield classes, 
\[\<X\> \w \<Y\> = \<X \w Y\>.\]

We always have $\<X\>\w\<Y\>\leq \<X\>\curlywedge\<Y\>$.

\begin{definition}  Define the following.  \label{def-dl-ba}
\begin{enumerate}
   \item Define $\DL =\left\{ \< X\>\in \BL \mbox{ with } \<X\> = \<X \w X\>\right\}$.
   \item A Bousfield class $\<X\>$ is called \textit{complemented} if there exists a class $\<X^c\>$ such that $\<X\> \w \< X^c\> =\<0\>$ and $\<X\> \vee \< X^c\> =\<\max\>$.  Call $\<X^c\>$ a \textit{complement} of $\<X\>$.
   \item Define $\BA$ to be the collection of Bousfield classes in $\DL$ that are complemented and have a complement in $\DL$.
\end{enumerate}
\end{definition}

When the category needs to be specified, we will write $\BL_\cC$, $\DL_\cC$, and $\BA_\cC$, or $\BL(\cC)$, etc.  In the case where $\cC=D(R)$ is the derived category of a ring, we will use the notation $\BL_R$, $\DL_R$, and $\BA_R$, or $\BL(R)$, etc.~instead.

The sub-poset $\DL\subseteq \BL$ is closed under arbitrary joins, and under the tensor operation, but not under meets; the meet in $\BL$ of two elements of $\DL$ may not be in $\DL$.  However, when we restrict to $\DL$, the meet is given by tensoring: if $\<X\>, \<Y\>,\<Z\> \in \DL$ have $\<Z\>\leq \<X\>$ and $\<Z\>\leq \<Y\>$, then $\<Z\> = \<Z\w Z\> \leq \<X \w Y\>$, so $\<X \w Y\>$ is the greatest lower bound.  A lattice is called \textit{distributive} if meets distribute across finite joins (equiv. if joins distribute across finite meets; see~\cite[I.6,~Thm.~9]{[Bir]}); it is a \textit{frame} if meets distribute across arbitrary joins.  Since the tensor product commutes with arbitrary coproducts, $\DL$ is a frame.

In general, a complemented class may have multiple complements.  When $\one\in\cC$, every complemented class is in $\DL$, because then $\<\max\>=\<\one\>$ and we have
\[ \<X\> = \<X\w\one\> = \<X\>\w (\<X\>\vee\<X^c\>) = (\<X\>\w \<X\>) \vee (\<X\>\w \<X^c\>) = \<X\w X\>.\]

Furthermore, if $\<X\>\in \BA$ then $\<X\>$ has a unique complement in $\DL$.  Indeed, if $\<X^c\>, \<\widetilde{X^c}\>\in \DL$ are two complements, then since the meet is given by tensoring, we have
\[ \<X^c\> = \<X^c\>\curlywedge(\<X\>\vee\<X^c\>) = \<X^c\>\w(\<X\>\vee\<\widetilde{X^c}\>) = \<X^c\>\w\<\widetilde{X^c}\>, \] and likewise $\<\widetilde{X^c}\> = \<X^c\>\w\<\widetilde{X^c}\>$.

One can check that $\BA$ is a sublattice of $\DL$ (i.e.~is closed under finite joins and meets), with $\<(X\vee Y)^c\> = \<X^c\> \w \<Y^c\>$ and $\<(X\w Y)^c\> = \<X^c\>\vee \<Y^c\>$.  In general, however, $\BA$ is not closed under infinite joins.  A \textit{Boolean algebra} is a distributive lattice in which every element is complemented; thus $\BA$ is a Boolean algebra, and this explains the notation.

We can use the tensor product to define another operation on Bousfield classes.

\begin{definition} For any Bousfield class $\<Z\>$ in the Bousfield lattice $\BL$, define the \textit{complementation operator}  $a(-)$ to be
\[ a\<Z\> = \bigvee_{\<Y\w Z\>=\0} \<Y\>.\]
\end{definition}

The complementation operator was first considered by Bousfield~\cite{[Bou79a]}, and later by Hovey and Palmieri~\cite{[HP]}, in the stable homotopy category.  Note that the definition requires knowing there is a set of Bousfield classes.  They prove the following properties of $a(-)$ in that context, but the proof is formal and applies in any well generated tensor triangulated category, or any well generated localizing tensor ideal of such a category.\\

\begin{lemma}~\cite[Lemma 2.3]{[HP]} \label{a-props}The complementation operator $a(-)$ has the following properties.
\begin{enumerate} 
   \item $\<E\>\leq a\<X\>$ if and only if $\<E\>\w \<X\>=\0$.
   \item $a(-)$ is order-reversing: $\<X\>\leq \<Y\>$ if and only if $a\<X\>\geq a\<Y\>$.
   \item $a^2\<X\>=\<X\>$.
\end{enumerate}
\end{lemma}

Note that we always have $\<X\>\w a\<X\>=\0$ and $\<X\> \vee a\<X\>\leq \<\max\>$.  If $\<X\>\in \DL$, then $a\<X\>$ is not necessarily in $\DL$.  If $\<X\>$ is complemented, with some complement $\<X^c\>$, then $\<X\>$ is also complemented by $a\<X\>$.  This is because, by the Lemma, $\<X\>\w \<X^c\>=\0$ implies $\<X^c\>\leq a\<X\>$, and thus $\<\max\>=\<X\>\vee\<X^c\> \leq \<X\>\vee a\<X\>$.  It follows that if $\one\in\cC$ and $\<X\>\in \BA$, then $a\<X\>$ is in $\DL$ and is the unique complement of $\<X\>$. 

We briefly mention a surprising but simple result using complementation, which we have been unable to find in the literature.  Call an object $X\in \cC$ \textit{square-zero} if $X$ is nonzero but $X\w X=0$.

\begin{proposition} \label{sq-free-1} Assume $\one\in\cC$. If there are no square-zero objects in $\cC$, then every object is complemented.
\end{proposition} 

\begin{proof}  Let $X\in \cC$ be arbitrary.  It suffices to show that $\<X\>\vee a\<X\>\geq \<\max\>$.  Suppose $Y$ has $\<Y\>\w \<X\>=\0$ and $\<Y\>\w a\<X\>=\0$.  Part (1) of Lemma~\ref{a-props} implies that $\<Y\>\leq \<X\>$, and from this we conclude that $Y\w Y=0$.  Our assumption forces $Y=0$ so $\<Y\>\w\<\max\>=\0$.  Thus $\<X\>$ is complemented by $a\<X\>$. \end{proof}

\begin{corollary} \label{sq-free-2} Assume $\one\in\cC$. If $\DL=\BL$, then $\BA=\DL=\BL$.
\end{corollary}

\begin{corollary} \label{sq-free-3} Assume $\one\in\cC$. There are no square-zero objects in $\cC$ if and only if \[\BA=\DL=\BL.\]
\end{corollary}

\subsection{Subcategories and quotient categories} Well generated categories behave well under taking subcategories and quotients.  A localizing subcategory $\cB\subseteq \cC$ is well generated if and only if $\cB = \loc(X)$ for some $X\in \cC$~\cite[Rmk.~2.2]{[IK]}.  Note that every Bousfield class is a localizing subcategory, and in fact a tensor ideal.

\begin{lemma} Every Bousfield class $\<Z\>\subseteq \cC$ is well generated.  Thus for all $Z\in \cC$, there exists an element $aZ\in \cC$ such that $\<Z\> = \loc(aZ)$.
\end{lemma}

\begin{proof} This follows from Proposition 2.1 in~\cite{[IK]}, since $\<Z\>$ is the kernel of the exact coproduct-preserving functor $F=(-\w Z)\colon \cC\to \cC$.
\end{proof}

\begin{lemma} For any $Z\in \cC$, we have $\<aZ\> = a\<Z\>$.
\end{lemma}

\begin{proof}  Because $aZ\in \<Z\>$, Lemma~\ref{a-props} implies that $\<aZ\>\leq a\<Z\>$.  If $\<Y\>$ has $\<Y\w Z\>=\0$, then $Y\in \<Z\>=\loc(aZ)$.  It follows that $\<Y\>\leq \<aZ\>$.  Therefore $a\<Z\>\leq \<aZ\>$ and equality holds.  \end{proof}

\begin{remark} \label{max-x} If $\cB\subseteq \cC$ is a well generated localizing tensor ideal, then $\cB=\loc(X)$ for some $X\in \cC$, and in this case $\<\max\>$ in $\BL(\cB)$ is $\<X\>$.  This is because $Y\in \loc(X)$ always implies $\<Y\>\leq \<X\>$.
\end{remark}

If $\cB\subseteq \cC$ is any localizing tensor ideal, we can form the Verdier quotient $\cC/\cB$.  This category has a tensor triangulated structure induced by that on $\cC$, such that the quotient functor $\pi\colon \cC\to \cC/\cB$ is exact, and $\pi(\one_\cC) = \one_{\cC/\cB}$.  If $\cB$ is well generated, then so is $\cC/\cB$, by~\cite[Cor.~4.4.3]{[nee-tricat]} or~\cite[Thm.~7.2.1]{[krause-loc-survey]}.  This also implies that $\cC/\cB$ has $\Hom$ sets.

\subsection{Examples} Next we survey several examples of categories and their Bousfield lattices.

\begin{example} Let $R$ be a commutative ring, or a graded-commutative ring.  Let $D(R)$ denote the unbounded derived category of right $R$-modules, or of  right graded $R$-modules (with degree-preserving maps).  If $R$ is graded, we think of objects in $D(R)$ as bi-graded; in either case we assume the differential decreases the chain degree by one.  Then $D(R)$ is a tensor triangulated category, with the product $A\w B = A\otimes_R^L B$ given by the left derived tensor product~\cite[\textsection 9.3]{[HPS]}.  The tensor unit is the module $R$ thought of as a complex concentrated in degree zero.  Furthermore, $D(R)=\loc(R)$, so $D(R)$ is compactly generated, hence well generated.  When $R$ is graded, this is meant in the multigraded sense discussed in~\cite[\textsection 1.3]{[HPS]}, and we follow the conventions of~\cite[\textsection 2]{[DP]}.  See also~\cite{[ivo-greg]}.  The Bousfield lattice of $D(R)$ is well-understood when $R$ is Noetherian; see the next example.  When $R$ is non-Noetherian our understanding of the Bousfield lattice is limited to several specific rings; see Example~\ref{example-lambda}. 
\end{example}

\begin{example} \label{example-strat} Iyengar and Krause~\cite{[IK]} investigate the Bousfield lattice of a compactly generated tensor triangulated category that is stratified by the action of a graded Noetherian ring $R$.  This general setting, developed in~\cites{[BIK_loc_cohom], [BIK_strat]}, building on~\cites{[Nee],[BCR],[HPS]}, includes the unbounded derived category of a commutative Noetherian ring; the stable module category StMod$(kG)$ of a finite group, where the characteristic of $k$ divides the order of the group, and then also the homotopy category $K(\mbox{Inj } kG)$ of complexes of injectives; and DG modules over a formal commutative DG algebra with a Noetherian cohomology ring.  They show that in such a category the Bousfield lattice is isomorphic to the lattice of subsets of the homogeneous prime spectrum of $R$, and  $\BA = \DL = \BL$.  In the case of a commutative Noetherian ring $R$, $D(R)$ is stratified by $R$, and so $\BL_R$ is isomorphic to the lattice of subsets of $\Spec R$.  The isomorphism is given in terms of support.
\end{example}

\begin{example} The ($p$-local) stable homotopy category $\S$ is a tensor triangulated category, with the product the smash product, and the unit the ($p$-local) sphere spectrum $S^0$.  Since $\S = \loc(S^0)$, this category is well generated.  Bousfield~\cite{[Bou79a]} showed that the class of every finite spectrum is in $\BA$, the class of every ring spectrum is in $\DL$, but for example the class of $H\bz$ is in $\DL$ but not in $\BA$.  He also showed that the Brown-Comenetz dual $IS^0$ of the sphere has $IS^0\w IS^0=0$, so $\DL\subsetneq \BL$.  Hovey and Palmieri~\cite{[HP]} study finer structure of the Bousfield lattice of this category.
\end{example}

\begin{example} \label{example-lambda} Fix a countable field $k$ and integers $n_i>1$, and consider the ring
\[\Ll= \frac{k[x_1,x_2,...]}{(x_1^{n_1}, x_2^{n_2},...)} ,\] with the $x_i$ graded so that $\Ll$ is graded-connected and finite-dimensional in each degree.  Let $D(\Ll)$ be the derived category of graded $\Ll$-modules; objects in $D(\Ll)$ are bigraded.  Neeman~\cite{[Nee-odd]} first considered such a ring (with $n_i=i$), showing the Bousfield lattice is large, although the homogeneous prime spectrum is trivial.  Dwyer and Palmieri~\cite{[DP]} examine the Bousfield lattice of $D(\Ll)$ in depth.  They show the Bousfield lattice has cardinality exactly $2^{2^{\aleph_0}}$.

Let $\IN = \Hom_k^*(\Ll,k)$ be the graded vector-space dual of $\Ll$.  This is a $\Ll$-module, and we consider it as an object of $D(\Ll)$ concentrated at chain degree zero.  The module $\IN$ plays an important role in~\cite{[DP]}.  One computation gives $\IN\w\IN = 0$, so $\DL_{\Ll}\subsetneq \BL_{\Ll}$.  This is relevant, because it implies that there is no Noetherian ring that stratifies $D(\Ll)$.

Furthermore, $\<\IN\>$ is a minimum nonzero Bousfield class: Corollary 7.3 in~\cite{[DP]} shows that for any non-zero $E$ in $D(\Ll)$, we have that $\<\IN\>\leq \<E\>$.  This implies that $\BA_{\Ll}$ is trivial, i.e.~the only complemented pair is $\<0\>$ and $\<\Ll\> = \<\one\>$ (see Prop.~\ref{no-min}).  
\end{example}

We mention one more difference among the Bousfield lattices in these examples.  One can easily check that every Bousfield class is a localizing subcategory.  Hovey and Palmieri~\cite[Conj.~9.1]{[HP]} conjecture that the converse holds in the stable homotopy category, but no progress has been made on this question.  In a category that is stratified by the action of a Noetherian ring, it is indeed the case that every localizing subcategory is a Bousfield class~\cite[Cor.~4.5]{[IK]}.  On the other hand, Greg Stevenson~\cite{[stevenson-counterexample]}, working in the unbounded derived category of a non-Noetherian ring (specifically any absolutely flat ring which is not semi-artinian), recently exhibited a localizing subcategory that is not a Bousfield class.

\subsection{Some (more) lattice theory.}

Here we recall some terminology and facts from lattice theory that we will need; our reference  is~\cite{[Bir]}.    A \textit{sub-poset} $\K$ of a lattice $\Lat$ is a subset of $\Lat$ along with the induced partial ordering.  A sub-poset $\K$ of a lattice $\Lat$ is a \textit{sublattice} if it is closed under finite joins and meets.

If $\K$ and $\Lat$ are lattices, a set map $F\colon\K\to \Lat$ is a \textit{join-morphism} if it is order-preserving (so $x\leq y$ implies $Fx\leq Fy$) and preserves binary joins.  A \textit{lattice morphism}  is a join-morphism that also preserves binary meets.  A \textit{lattice isomorphism} is a lattice morphism that is a set bijection and has an order-preserving inverse.

We do not assume that a  join-morphism preserves minimum or maximum elements.  Nor do we assume that a join-morphism between Bousfield lattices will commute with the tensor product operation $\<X\>\w \<Y\>$.

Note that if $F$ is a bijection with inverse $G$, and both $F$ and $G$ are join-morphisms that preserve arbitrary joins, then they preserve binary meets so $F$ is a lattice isomorphism.

Any poset can be thought of as a category, where $x\leq y$ if and only if there is a (unique) morphism from $x$ to $y$.  Joins are colimits and meets are limits.  Then a complete lattice corresponds to a category that is complete and cocomplete in the categorical sense.  

\begin{definition} For any element $a$ in a lattice $\Lat$, define $a\da = \{x\in \Lat\;|\; x\leq a\}$ and $a\ua = \{x\in \Lat\;|\; x\geq a\}$.  Note that these are both sublattices of $\Lat$.
\end{definition}

\begin{definition} A nonempty subset $J$ of a complete lattice $\Lat$ is an \textit{ideal} if it is closed under finite joins, and $a\in J$ and $x\in \Lat$ with $x\leq a$ implies $x\in J$.  An ideal is \textit{complete} if it is closed under arbitrary joins.  Note that $a\da$ is an ideal, for all $a\in \Lat$.  An ideal $J$ is \textit{principal} if $J=a\da$ for some $a\in\Lat$.  Note that an ideal $J$ is principal if and only if it is complete.   \end{definition}

\begin{definition} \label{defn-quotient} Given a principal ideal $J$ of a complete lattice $\Lat$, and $a,b\in \Lat$, we say $a\equiv b \modJ$, and write $[a]=[b]$,  if $a\vee c = b\vee c$ for some $c\in J$.  The equivalence classes under this equivalence relation, with the ordering, join, and meet induced by $\Lat$, form a complete lattice $\Lat/J$, called the \textit{quotient lattice}.  The quotient map $\Lat\to \Lat/J$ sending $x\mapsto [x]$ is a lattice epimorphism. \end{definition}

It is not hard to show that if $J=a\da$ is a principal ideal in a complete lattice $\Lat$, then $[x]=[y]$ in $\Lat/J$ if and only if $x\vee a = y\vee a$.  Every quotient of a complete lattice by a principal ideal is isomorphic to a sublattice: for all $a\in \Lat$, there is an isomorphism of lattices $\Lat/a\da \stackrel{\sim}{\to} a\ua$, given by $[x]\mapsto x\vee a$.

\begin{definition} Given lattices $\K$ and $\Lat$, the product lattice is defined as the set product $\K\times\Lat$, with $(a,b)\leq(c,d)$ precisely when $a\leq c$ and $b\leq d$, and joins and meets  defined termwise.  One can check that, for example, $0\times \Lat$ is a principal ideal in $\K\times \Lat$, and there is a lattice isomorphism $(\K\times\Lat)/(0\times \Lat)\cong \K$.
\end{definition}


\section{Lattices and quotients} \label{section-landq}

In this section we give some results comparing the quotient of a Bousfield lattice to the Bousfield lattice of a quotient.  Again, let $\cC$ be a well generated tensor triangulated category, or a well generated localizing tensor ideal of such a category.  Let $Z$ be an element of $\cC$, and consider the Verdier quotient $\cC/\<Z\>$ and quotient functor $\pi\colon \cC\to \cC/\<Z\>$.

\begin{lemma} The functor $\pi$ induces an onto join-morphism of lattices that preserves arbitrary joins,
\[ \pi\colon \BL(\cC)\to \BL\left(\cC/\<Z\>\right), \mbox{ where } \<X\> \mapsto \<\pi X\>. \]
\end{lemma}

\begin{proof}  We will show that if $X, Y\in\cC$ have $\<X\>\leq \<Y\>$ in $\BL_\cC$, then $\<\pi X\>\leq \<\pi Y\>$ in $\BL\left(\cC/\<Z\>\right)$.  This will show that $\pi$ is order-preserving, and by symmetry will also show that $\pi$ is well-defined.  Take $W\in \cC/\<Z\>$ with $W\w \pi Y=0$.  Take $\widetilde{W}$ in $\cC$ so $\pi\widetilde{W} = W$.  The tensor structure on $\cC/\<Z\>$ is such that  $\pi(\widetilde{W}\w Y)=\pi\widetilde{W}\w\pi Y=0$, so we have $\widetilde{W}\w Y\in \<Z\>$, i.e.~$\widetilde{W}\w Y\w Z=0.$  Then $(\widetilde{W}\w Z) \w X = 0$, by hypothesis, so $\pi(\widetilde{W}\w X)=0$.  This shows that $W\w \pi X=0$. 

Since $\pi$ commutes with arbitrary coproducts, it commutes with arbitrary joins.\end{proof}

\begin{proposition} \label{epim} Assume $\one\in\cC.$ Let $\<Z\>$ be any Bousfield class in $\BL_\cC$.  Then $\pi$ induces an onto join-morphism of lattices that preserves arbitrary joins,
\[  \overline{\pi}\colon \BL_\cC/(a\<Z\>)\da\lra \BL\left(\cC/\<Z\>\right) , \] such that if $\overline{\pi}[\<X\>] = \0$, then $[\<X\>]=[\0]$.
\end{proposition}

\begin{proof}  First we show that $\overline{\pi}$ is order-preserving.  Suppose $\<X\>\leq \<Y\>$ in $\BL_\cC/(a\<Z\>)\da$; this is equivalent to assuming $\<X\>\vee a\<Z\> \leq \<Y\>\vee a\<Z\>$.  We want to show that $\<\pi X\>\leq \<\pi Y\>$ in $\BL(\cC/\<Z\>)$.  Take $W\in \cC/\<Z\>$ with $W\w \pi Y=0$, and let $\widetilde{W} \in \cC$ be such that $\pi \widetilde{W} = W$.  Then $0=\pi\widetilde{W}\w \pi Y = \pi (\widetilde{W} \w Y)$, so $\widetilde{W}\w Y\in \<Z\>$ and $\widetilde{W}\w Y \w Z=\widetilde{W} \w Z \w Y=0$.

Since $\<Z\>\w a\<Z\> =\0$, we also have $\<\widetilde{W} \w Z\>\w a\<Z\>=\0$.  Our assumption then implies that $(\widetilde{W}\w Z) \w X=0$.  Therefore $\widetilde{W}\w X\in \<Z\>$, which says that $0=\pi(\widetilde{W}\w X) = W\w \pi X$.

Thus $\overline{\pi}$ is order-preserving, and hence well-defined.  It is clearly onto and preserves joins.

Now suppose $\overline{\pi}[\<X\>]=\0$.  Then $\<\pi X\> = \0$, so $X\in \<Z\>$, i.e.~$X\w Z=0$.  Using Lemma~\ref{a-props} this implies that $\<X\>\leq a\<Z\>$, so $[\<X\>]=[\0]$.   \end{proof}

To be a lattice isomorphism, $\overline{\pi}$ must have an order-preserving inverse.  In the remainder of this section, we will give examples of when this does and does not happen.

\begin{corollary} \label{compl} Let $\<Z\>$ and $\<Z^c\>$ be a pair of complemented classes in $\BL_\cC$, and assume that $\one\in\cC$.  Then $\pi$ induces a lattice isomorphism
\[ \overline{\pi}\colon \BL_\cC/\<Z^c\>\da\lra \BL\left(\cC/\<Z\>\right). \]
\end{corollary}

\begin{proof}  Recall that $\one\in\cC$ implies that complements are unique and $\<Z^c\>=a\<Z\>$.  We claim that $\psi(\<\pi X\>) = [\<X\>]$ is a well-defined, order-preserving inverse to $\overline{\pi}$.  Thus we wish to show that if $\<\pi X\> \leq \<\pi Y\>$, then $\<X\>\vee \<Z^c\> \leq \<Y\>\vee \<Z^c\>$.  By symmetry, this will also show that $\psi$ is well-defined; by inspection, then, it is an inverse to $\overline{\pi}$.

Take $W\in \cC$ with $W\in \<Y\>\vee \<Z^c\>$.  Then $W\w Y=0$, so $W\w Y\w Z=0$, i.e.~$(W\w Y)\in \<Z\>$.  This says that $\pi(W\w Y)=0$ in $\cC/\<Z\>$, so $\pi W\w \pi Y=0$.  By hypothesis, this means $\pi W\w \pi X=0$, which working backwards implies that $W\w X\w Z=0$.

On the other hand, we also know that $W\w Z^c=0$, so $W\w X\w Z^c=0$.  Therefore $(W\w X) \in \<Z\>\vee \<Z^c\> =\<\max\>= \<\one\>$, so $W\w X=0$.  Thus $W\in \<X\vee Z^c\>$ as desired. \end{proof}

For example, if $L\colon\cC\to\cC$ is a smashing localization functor with colocalization $C$, and $\one\in \cC$, then $\<L\one\>$ and $\<C\one\>$ are a complemented pair.  This result relates the Bousfield lattice of $\cC$ to the Bousfield lattice of the $L$-local category, which is equivalent to $\cC/\<L\one\>$.  See also~\cite[Prop.~6.12]{[IK]}.

\begin{corollary} Suppose $\BL_\cC = \DL_\cC$.  Then for every Bousfield class $\<Z\>$, the functor $\pi$ induces a lattice isomorphism
\[ \overline{\pi}\colon \BL_\cC/(a\<Z\>)\da\lra \BL\left(\cC/\<Z\>\right).\]
\end{corollary}

\begin{proof} This follows immediately from Corollaries~\ref{sq-free-2} and~\ref{compl} if $\one\in\cC$, but we will prove it more generally. As in the last proof, we claim that $\psi(\<\pi X\>) = [\<X\>]$ is a well-defined, order-preserving inverse to $\overline{\pi}$.  Suppose $\<\pi X\>\leq \<\pi Y\>$; it suffices to show that $\<X\>\vee a\<Z\>\leq \<Y\>\vee a\<Z\>$.  

Take $W\in \cC$ with $W\in \<Y\>\vee a\<Z\>$.  As in the last proof, $W \w Y=0$ implies $W\w X\w Z=0$.  Then $\<W\>\w a\<Z\>=\0$ implies $\<W\>\leq \<Z\>$, by Lemma~\ref{a-props}.  Therefore $(W\w X)\in \<W\>$.  Since $\BL_\cC=\DL_\cC$, we have $X\in \<W\w W\>=\<W\>$, so $W\w X=0$ and this concludes the proof. \end{proof}

The previous two corollaries apply when $\cC$ is a stratified category, as discussed in Example~\ref{example-strat}.  The next result, however, shows that $\overline{\pi}$ fails to be an isomorphism in both the stable homotopy category and the category $D(\Ll)$.

\begin{proposition} \label{not-isomo} Assume $\one\in \cC$.  Suppose $\<Z\>$ in $\BL_\cC$ is an element of $\DL_\cC\backslash \BA_\cC$.  Then the map induced by $\pi$ in Proposition~\ref{epim} is not an isomorphism.  This happens in both the stable homotopy category and $D(\Ll)$.
\end{proposition}

\begin{proof}  Since $\<Z\>\notin \BA_\cC$, we know $\<Z\>\vee a\<Z\> < \<\one\>$.  We will show that $\<\pi Z\>= \<\pi \one\>$ in $\BL(\cC/\<Z\>)$, but $[\<Z\>] \neq [\<\one\>]$ in $\BL_\cC/(a\<Z\>)\da$.

Since $\pi$ is order-preserving, we know that $\<\pi Z\>\leq \<\pi \one\>$.  We must show $\<\pi Z\>\geq \<\pi \one\>$.  Suppose $W\in \cC/\<Z\>$ has $W\w \pi Z=0$.  Choose $\widetilde{W}\in \cC$ such that $\pi \widetilde{W} = W$.  Then $\pi(\widetilde{W}\w Z)=0$, so $\widetilde{W}\w Z\w Z=0$.  This says $\widetilde{W}\in \<Z\w Z\>$, and by hypothesis $\<Z\w Z\>=\<Z\>$, so $\widetilde{W}\w Z=0$.  Therefore $\widetilde{W}=\widetilde{W}\w\one\in \<Z\>$ and $\pi (\widetilde{W}\w\one) = W\w\pi\one = 0$ in $\cC/\<Z\>$.  This shows $\<\pi\one\> = \<\pi Z\>$.

By assumption, $\<Z\>$ is such that $\<Z\>\vee a\<Z\> < \<\one\>$.  But $\<\one\> = \<\one\> \vee a\<Z\>$, so $\<Z\>\vee a\<Z\> < \<\one\>\vee a\<Z\>$ and thus $[\<Z\>] <[\<\one\>] $ in $\BL_\cC/(a\<Z\>)\da$. 

In the $p$-local stable homotopy category, we can take $\<Z\> = \<H\Fp\>$.  The spectrum $H\Fp$ is a ring spectrum, and Bousfield~\cite{[Bou79a]} shows that the Bousfield class of any ring spectrum is in the distributive lattice.  Let $IS^0$ be the Brown-Comenetz dual of the sphere.  Then Lemma 7.1 in~\cite{[HP]} shows that $IS^0\w H\Fp=0$, and $\<IS^0\>\leq \<H\Fp\>$ so $\<IS^0\>\w a\<H\Fp\>=\0$ by Lemma~\ref{a-props} above.  This shows that $IS^0\in \<H\Fp\>\vee a\<H\Fp\>$, and hence $\<H\Fp\>\vee a\<H\Fp\> < \<S^0\>$.

In the category $D(\Ll)$ of Example~\ref{example-lambda}, we can take $\<Z\> = \<k\>$.  The class $\<k\>$ is in $\DL_\Ll$ because $k$ is a ring object.  The dual $\IN$ of $\Ll$ has $\IN\w k=0$~\cite[Cor.~4.12]{[DP]} and $\<\IN\>\leq \<k\>$~\cite[Lemma 4.8]{[DP]}, so $\<\IN\>\w a\<k\>=\0$.  Thus we have that $\IN\in \<k\>\vee a\<k\> < \<\Ll\>$.   \end{proof}


\section{Ring maps and the Bousfield lattice} \label{section-rmbl}

\subsection{Ring maps and derived categories}

In this section, we'll establish basic facts about morphisms of Bousfield lattices induced by ring maps, laying the groundwork for the results in Section~\ref{section-nnr}.

\noindent \textbf{WARNING: The results of this section hold in an ungraded or a graded setting, and we will be ambiguous with notation.}  Thus let $f\colon R\to S$ be either a ring homomorphism between two commutative rings, or a graded ring homomorphism between two graded-commutative rings.  Let $\Mod$-$R$ denote either the category of right $R$-modules, or the category of right graded $R$-modules (with degree-preserving maps).  Let $D(R)$ denote either the unbounded derived category of $R$-modules, or the unbounded derived category of graded $R$-modules; in the latter case, the objects of $D(R)$ are bigraded, and we follow the conventions in~\cite[\textsection 2]{[DP]}.\\

In either case we will use the standard model structure on the category $Ch(R) = Ch(\Mod$-$R)$ of unbounded chain complexes.  The weak equivalences are quasi-isomorphisms, the fibrations are dimensionwise surjections, and the cofibrations are dimensionwise injections with cofibrant cokernels.  The cofibrant objects are the complexes that can be written as an increasing union of subcomplexes such that the associated quotients are complexes of projectives with zero differentials.  Every object is fibrant.  See~\cite[\textsection 2.3]{[hovey-model]} or~\cite[\textsection 9.3]{[HPS]} for more details.  

A ring map $f\colon R\to S$ induces a functor on module categories $f_*\colon \Mod$-$R \to \Mod$-$S$, via extension of scalars, where $f_*(M) = M\otimes_RS$.  This induces a functor $f_*\colon Ch(R)\to Ch(S)$ on chain complexes.  The forgetful functor $f^*\colon\Mod$-$S\to \Mod$-$R$ induces a functor $f^*\colon Ch(S)\to Ch(R)$, and $f_*$ and $f^*$ are adjoints.

\begin{definition} Let $\f $ be the left derived functor $\f  = Lf_* = L(-\otimes_R S)\colon D(R)\to D(S)$.  Let $\i = Rf^*\colon D(S)\to D(R)$ be the right derived functor of the forgetful functor.  \end{definition}

\begin{lemma} The derived functors $\f$ and $\i$ exist and form a Quillen adjoint pair; $\f$ is the left adjoint and $\i$ is the right adjoint.
\end{lemma}

\begin{proof}  Since $f_*$ is left adjoint to $f^*$, by~\cite[Rmk.~9.8]{[DS]} it suffices to show that $f^*$ preserves fibrations and trivial fibrations.  Since fibrations are degreewise surjections, this is immediate.
\end{proof}

The functor $\f $ is exact (i.e.~sends exact triangles to exact triangles), has $\f (R)=S$, and $\f (X\w Y) = \f X\w \f Y$ (see~\cite[Thm.~9.3.1]{[HPS]} and note that they consider both the ungraded and graded settings).  Since it is a left adjoint, it commutes with coproducts.  Since every object is fibrant, we have $\i(X)=f^*(X)$ for all $X$, so $\i$ is exact and commutes with coproducts and products.  

\begin{remark} Take $z\in R_0$, and consider the morphism $R\stackrel{z}{\to} R$ in $D(R)$.  Applying $\f $ to this, we get
\[ \left(\f (R)\stackrel{\f (z)}{\lra}\f (R)\right) = \left(R\otimes_RS\stackrel{z\otimes1}{\lra} R\otimes_RS\right) \]
\[  = \left(R\otimes_RS\stackrel{1\otimes f(z)}{\lra} R\otimes_RS\right) = \left(S\stackrel{f(z)}{\lra} S\right).\\ \\\]
\label{Gz}\end{remark}

The following lemma, called the projection formula and proved in~\cite{[Wei]} for bounded-below complexes, will be used frequently. 

\begin{lemma}{(Projection Formula)} \label{proj}  For all objects $A$ in $D(R)$ and $B$ in $D(S)$, we have
\[ \i(\f A\w B) = A\w \i B. \]
\end{lemma}

\begin{proof} Recall that we can compute the derived tensor product $-\w-$ by taking a cofibrant replacement in either factor.  Let $QX$ represent a choice of cofibrant replacement for a complex $X$.

Since every object is fibrant, we have
\[ \i(\f A\w B) = f^* (\f A\w B) = f^*(Q(\f A)\otimes_S B).\]

To compute $\f A$ we use a cofibrant replacement $QA$ of $A$.  Since $f_*$ is left Quillen, it preserves cofibrant objects.  Thus $Q(\f A) = Q(f_*(QA)) = f_*(QA)$.

At the module level, for $M\in \Mod$-$R$ and $N \in \Mod$-$S$, we have
\[  f^*\left(f_*(M)\otimes_S N\right) = f^*\left((M\otimes_RS)\otimes_SN\right)=M\otimes_Rf^*(N), \] and this extends to the level of chain complexes, to give
\[ \i(\f A\w B) = f^*(f_*(QA)\otimes_S B) = QA\otimes_R f^*(B) = A\w \i B. \]  
\end{proof}

\begin{corollary} \label{proj_bare} For all objects $A$ in $D(R)$ and $B$ in $D(S)$,
\[ \f A\w B = 0 \mbox{ if and only if }  A\w \i B=0. \]
\end{corollary}

\subsection{Induced maps on Bousfield lattices}

Here we show that the functors $\f $ and $\i$ induce maps between the Bousfield lattices of $D(R)$ and $D(S)$.  If we consider a Bousfield class $\<X\>$ as the localizing subcategory of $X$-acyclics, then we can map this to $\f (\<X\>)$ as a subcollection in $D(S)$; however, in general $\f (\<X\>)$ will not be triangulated.   Instead we make the following definitions.

\begin{definition} Define a map $\f \colon\BL_{R}\to \BL_{S}$ by $\<X\>\mapsto \<\f X\>$.  Also,  define a map $\i\colon\BL_{S}\to \BL_{R}$ by $\<X\>\mapsto \<\i X\>$.  For the rest of this document, $\f \<X\>$ and $\i\<X\>$ will mean $\<\f X\>$ and $\<\i X\>$.\label{def-gxix} \end{definition}

\begin{proposition} \label{G_i_def}  Both $\f $ and $\i$ induce join-morphisms on Bousfield lattices that commute with arbitrary joins.
\end{proposition}

\begin{proof}  First we show that $\<Y\>\leq\<X\>$ implies $\<\i Y\>\leq\<\i X\>$.  Suppose $\<Y\>\leq\<X\>$ and $W\w \i X=0$.  Then Corollary \ref{proj_bare} implies  $\f W\w X=0$.  Thus $\f W\w Y=0$, and $W\w \i Y=0$.

This implies that if $\<Y\> = \<X\>$, then $\<\i Y\> = \<\i X\>$, so $\i$ is well-defined and order-preserving.

Now suppose $\<Y\>\leq\<X\>$ and $\f X\w W=0$.  Then from Corollary \ref{proj_bare}, $X\w \i W =0$, so $Y\w \i W=0$, which implies $\f Y\w W=0$.  Therefore $\f$ is order-preserving and well-defined.  Both $\f $ and $\i$ commute with coproducts on the object level, hence with arbitrary joins at the level of Bousfield classes.  \end{proof}

Note that $\f$ commutes with the tensor operation, $\f(\<X\w Y\>) = \<\f X\> \w \<\f Y\>$, but in general $\i$ does not.  See Lemma~\ref{GiX_equiv} however.

Recall from Section~\ref{section-BLbackground} that since $R\in D(R)$, $\<\max\>=\<R\>$ in $\BL_R$, and every complemented class in $\BL_R$ is in $\DL_R$.  Furthermore, complements are unique, and are given by the complementation operator $a(-)$.

\begin{lemma} The functor $\f $ maps $\DL_{R}$ into $\DL_{S}$, and $\BA_{R}$ into $\BA_{S}$.  If $\<X\>$ in $\BA_{R}$ has complement $\<\Xc \>$, then $\<\f X\>$ has complement $ \<\f (\Xc) \>$.
\label{f_DLBA}\end{lemma}

\begin{proof}  If $\<Y\> = \<Y\w Y\>$, then $\<\f Y\> = \<\f Y\w \f Y\>$.

If $\<X\>$ has $\<X\>\vee\<\Xc \>=\<R\>$ and $\<X\>\w\<\Xc \>=\<0\>$, then $\<\f X\>\vee \<\f (\Xc) \>=\<\f R\>=\<S\>$ and $\<\f X\>\w \<\f (\Xc) \>=\<0\>$, so $\<(\f X)^c\> = \<\f (\Xc) \>$.
\end{proof}

We will strengthen and extend this lemma in the next subsection, under additional hypotheses.  Next we describe a useful quotient of $\BL_{R}$.

\begin{definition} \label{def-jm} Fix $J_f$ to be the image of $Ker \f $ in $\BL_{R}$, in other words $J_f = \left\{\<X\>\;|\; \f \<X\>=\<0\>\right\}.$ Also define 
\[ \<M_f\> = \bigvee_{\<Y\>\in J_f}\<Y\>. \]
\end{definition}

\begin{proposition} \label{BLJ}  The subposet $J_f$ is a principal ideal in $\BL_{R}$ with $J_f=\<M_f\>\da$, and $\f $ induces a lattice join-morphism that preserves arbitrary joins,
\[ \oG \colon \BL_{R}/J_f\to \BL_{S}, \] where $\oG\<X\> = \<\f X\>$.
\end{proposition}

\begin{proof}  Suppose $\<Y\>\leq\<X\>$ and $\<\f X\>=\<0\>$.  Then $\<\f Y\>\leq\<\f X\>$, so $\<\f Y\>=\<0\>$ and $J_f$ is a lattice ideal.  Every $\<X\>$ in $J_f$ has $\<X\>\leq \<M_f\>$, so $J_f\subseteq \<M_f\>\da$.  And since $\f\<M_f\>=0$, if $\<X\>\leq \<M_f\>$ then $\<X\>\in J_f$.  Therefore $J_f=\<M_f\>\da$ is principal.

To get an induced map on the quotient lattice, we need to know that if $[\<X\>]= [\<Y\>]$, then $\f \<X\>=\f \<Y\>$.  Since $\<X\>$ and $\<Y\>$ are equivalent if and only if $\<X\>\vee \<M_f\> = \<Y\>\vee\<M_f\>$, and $\<\f M_f\>=\<0\>$, we get
\[ \<\f X\> = \<\f X\>\vee \<\f M_f\> = \f (\<X\>\vee \<M_f\>) = \<\f Y\>\vee\<\f M_f\> = \<\f Y\>. \] 

Thus $\overline{\f}$ is well-defined.  It is order and join-preserving since $\f$ is.
\end{proof}

\begin{remark} \label{mf-ais} Note that for any object $X$ in $D(R)$, the projection formula implies $\f X=0$ if and only if $X\w \i S=0$, which is true if and only if $\<X\>\leq \<M_f\>$.  Therefore by definition $\<M_f\> = a\<\i S\>$, and thus also $\<\i S\> = a\<M_f\>$.  \end{remark}

\begin{proposition} \label{rs-lat-isomo} If $\<\i S\>\vee \<M_f\> = \<R\>$, then the quotient functor $\pi\colon D(R)\to D(R)/\<\i S\>$ induces a lattice isomorphism
\[ \overline{\pi}\colon \BL_R/J_f \lra \BL(D(R)/\<\i S\>). \]
\end{proposition}

\begin{proof} This follows easily from Remark~\ref{mf-ais} and Corollary~\ref{compl}, and the fact that $J_f = \<M_f\>\da$. \end{proof}

We don't have a general criteria for when to expect $\<\i S\>\vee \<M_f\> = \<R\>$ to hold.  It holds when $\BL_R = \DL_R$, thanks to Corollary~\ref{sq-free-2}.  This is the case if $R$ is Noetherian, for example.  On the other hand, consider the ring $\Ll$ from Example~\ref{example-lambda}, and let $f\colon\Ll\to k$ be projection onto the degree zero piece.  Then $\<\i S\>\vee \<M_f\> = \<k\>\vee a\<k\>$, and as shown in the proof of Proposition~\ref{not-isomo}, this is strictly less than $\<\Ll\>$.

\subsection{Maps $f\colon R\to S$ satisfying $\f \i\<X\>=\<X\>$ for all $X$} \label{rmbl-fiX-X}

In this subsection we assume the map $f\colon R\to S$ satisfies $\f \i\<X\>=\<X\>$ for all $X$.   In Section~\ref{section-nnr} we show that this condition holds for the specific map of non-Noetherian rings $\g\colon \G\to\L$.

\begin{lemma} The following are equivalent:
\begin{enumerate}
	\item $\f \i\<X\>=\<X\>$ for all $\<X\>,$
	\item $\i W\w \i X = 0 \mbox{ if and only if } \i(W\w X)=0,$
	\item $\i\<Y\w X\> = \<\i Y\>\w \<\i X\>$.
\end{enumerate}
\label{GiX_equiv}\end{lemma}

\begin{proof}
For $(1)\Leftrightarrow (2)$, note that $W\w \f \i X=0$ iff $\i W\w \i X=0$, and  $W\w X=0$ iff $\i(W\w X)=0$.

For $(1)\Leftrightarrow (3)$, note that $W\w \i(Y\w X)=0$ iff $\f W\w (Y\w X) =0$ iff $(\f W\w Y)\w X=0$, and $W\w \i X\w \i Y=0$ iff $(\f W\w \f \i X)\w Y=0$ iff $(\f W\w Y)\w (\f \i X)=0$.
\end{proof}

This is a good setting in which to consider the behavior of the sub-posets $\BA$ and $\DL$ under $\f$ and $\i$.

\begin{lemma} Suppose $\f \i\<X\>=\<X\>$ for all $\<X\>$.  Then the map $\i$ sends $\BA_S$ into $\BA_R$ if and only if $\<\i S\>\vee\<M_f\> = \<R\>$.  If this is the case, and $\<X\>\in \BA_{S}$ has complement $\<\Xc \>$, then $\<\i X\>\in \BA_{R}$ has complement $\<\i (\Xc) \>\vee \<M_f\>$.
\label{i_BA}\end{lemma}

\begin{proof}  If $\i$ injects $\BA_S$ into $\BA_R$, then since $\<0\>$ and $\<S\>$ are a complemented pair in $\BA_S$, the class $\<\i S\>$ is complemented in $\BA_R$.  Its complement must be $a\<\i S\>$, which is $\<M_f\>$ by Remark~\ref{mf-ais}.

For the converse, suppose that  $\<\i S\>\vee\<M_f\> = \<R\>$.  Now suppose $\<X\>\in \BA_{D(S)}$, so $\<X\>\vee\<\Xc \>=\<S\>$ and $\<X\>\w\<\Xc \>=\<0\>$.  This implies $\<\i X\>\vee\<\i (\Xc) \>=\<\i S\>$ and $\<\i X\>\w\<\i (\Xc) \>=\<0\>$, using Lemma~\ref{GiX_equiv}. 

We calculate that
\[ \<\i X\>\vee(\<\i (\Xc) \>\vee\<M_f\>) = \<\i S\>\vee \<M_f\> =  \<R\>. \]

Also, we have
\[ \<\i X\>\w(\<\i (\Xc) \>\vee\<M_f\>) = (\<\i X\>\w\<\i (\Xc) \>) \vee (\<\i X\>\w\<M_f\>) \]
\[ = \<0\>\vee (\<\i X\>\w\<M_f\>)=\<0\>. \]

The last equality follows from the fact that $\i$ is order preserving and $\<X\>\leq \<S\>$ for all $\<X\>$, so $\<\i S\> \w \<M_f\>=\0$ implies $\<\i X\>\w \<M_f\>=\0$ for all $X$ in $D(S)$.

This shows that the complement of $\<\i X\>$ is $\<\i (\Xc) \>\vee \<M_f\>$.
\end{proof}

\begin{proposition} \label{G_i_sl} Suppose $\f \i\<X\>=\<X\>$ for all $\<X\>$.  The following hold.
\begin{enumerate}
    \item The map $\f $ sends $\DL_{R}$ onto $\DL_{S}$, and the map $\i$ injects $\DL_{S}$ into $\DL_{R}$.
    \item The map $\f$ sends $\BA_{R}$ onto $\BA_{S}$, and if $\<\i S\>\vee\<M_f\>=\<R\>$ then  $\i$ injects $\BA_{S}$ into $\BA_{R}$.
\end{enumerate}    
\end{proposition}

\begin{proof} Lemma \ref{GiX_equiv} implies that if $\<Y\> = \<Y\w Y\>$, then $\<\i Y\> = \<\i Y\w \i Y\>$, so $\i$ sends $\DL_S$ to $\DL_R$ and is injective by hypothesis.  The rest follows from Proposition \ref{f_DLBA}, Lemma \ref{i_BA}, and the fact that $\f$ is surjective and  $\i$ is injective.
\end{proof}

This is also a good setting in which to consider poset adjoints.  As a poset map, because $\i$ preserves joins on $\BL_{S}$, it has a poset map right adjoint $r\colon \BL_{R}\to\BL_{S}$,  see~\cite[Lemma 3.5]{[HP]}.  We know
\[ r\<Y\> = \bigvee \left\{\<X\> \; | \; \i\<X\>\leq \<Y\>\right\}, \mbox{ and } \]
\[ \i\<X\>\leq \<Y\> \mbox{ if and only if } \<X\>\leq r\<Y\>. \]

\begin{proposition} If $\f \i\<X\>=\<X\>$ for all $\<X\>$, then $\f \<X\> = r\<X\>$ for all $\<X\>$, so
\[ \<\i X\>\leq \<Y\> \mbox{ if and only if } \<X\>\leq \<\f Y\>. \]
\end{proposition}

\begin{proof} First suppose that $\<\i X\>\leq \<Y\>$ and $W\w \f Y=0$ for some $W$.  Then Corollary \ref{proj_bare} implies $\i W\w Y=0$, so $\i W\w \i X=0$.  It follows from Lemma \ref{GiX_equiv} that  $\i (W\w X)=0$, so $W\w X=0$.

For the other direction, if  $\<X\>\leq \<\f Y\>$, then $\<\i X\>\leq \<\i\f Y\>$.  In general we always have $\<\i\f Y\>\leq \<Y\>$.  Indeed, an object $W$ has $W\w \i\f Y=0$ iff $\f W\w \f Y=0$ iff $\f (W\w Y)=0$, so $W\w Y=0$ implies $W\w (\i \f  Y)=0$.

This immediately implies that  $\bigvee \left\{\<X\> \; | \; \i\<X\>\leq \<Y\>\right\} \leq \f\<Y\>$.  And the fact that  $\i\<\f Y\>\leq \<Y\>$ gives $\<\f Y\> \leq \bigvee \left\{\<X\> \; | \; \i\<X\>\leq \<Y\>\right\}$.
\end{proof}

The $\BL$ operation $\f $ also preserves arbitrary joins, so has a poset map right adjoint.  On the object level, we know that $\i $ is right adjoint to $\f$, and so it is natural to ask if $\i$ is the right poset adjoint of $\f$.

\begin{proposition} Assume $\f \i\<X\>=\<X\>$ for all $X$.  Then on the level of Bousfield classes, we have
\[ \<\f X\>\leq \<Y\> \Leftarrow \<X\>\leq \<\i Y\>, \] but the forward direction need not hold.  
\end{proposition}

\begin{proof}  First suppose $\<X\>\leq \<\i Y\>$ and $W\w Y=0$.  Then $\i(W\w Y)=0$, which using Lemma \ref{GiX_equiv} means $\i W\w \i Y=0$, so $\i W\w X=0$, and $W\w \f X=0$.

On the other hand, suppose $\<\f X\>\leq \<Y\>$ and $W\w \i Y=0$.  Then $\f W\w Y=0$, $\f W\w \f X=0$, and $\f (W\w X)=0$.  At the $\BL$ level, this does not necessarily mean $W\w X=0$.  (Take, for example, $Y=0$, $W=R$, and $X$ any object such that $\f X=0$.)  
\end{proof}

We end this section with another lattice isomorphism.

\begin{theorem} \label{latt-o-q} Suppose $\f \i\<X\>=\<X\>$ for all $\<X\>$.  There is a lattice isomorphism
\[ \phi\colon \BL\left( D(R)/\<\i S\>\right) \to \BL_S, \] given by $\phi\<X\> = \<\f\tilde{X}\>$, where $\pi\tilde{X}=X$.
\end{theorem}

\begin{proof} Recall that $\pi\colon D(R)\to D(R)/\<\i S\>$ is the canonical projection.  First we suppose that $\<X\>,\<Y\>$ in $\DRiS$ have $\<X\>\leq \<Y\>$, and will show that then $\phi\<X\>\leq \phi\<Y\>$.  Fix a choice of $\tilde{X}$ and $\tilde{Y}$ such that $\pi\tilde{X}=X$ and $\pi\tilde{Y}=Y$.

Take $W\in D(S)$ such that $W\w \f\tilde{Y}=0$.  We wish to show that $W\w\f\tilde{X}=0$.  Corollary~\ref{proj_bare} implies that $\i W\w \tilde{Y}=0$.  So
\[ 0=\pi(\i W\w\tilde{Y}) = \pi\i W\w \pi\tilde{Y} = \pi \i W\w Y. \]

By hypothesis, this implies that $\pi\i W\w X=0$.  Thus $\pi(\i W\w \tilde{X})=0$, and $\i W\w \tilde{X}\in \<\i S\>$, so $\i W\w \tilde{X}\w \i S=0$.  Again the projection formula implies that $\f(\i W\w \tilde{X})=0$, so $\f\i W\w \f\tilde{X}=0$.  Since we're assuming $\<\f\i W\> = \<W\>$, we conclude that $W\w\f\tilde{X}=0$, as desired.

This shows that $\phi$ is order-preserving.  By symmetry, it also shows that $\phi$ is well-defined, independent of choice of representative or preimage.  

The map $\phi$ is surjective by assumption: given $\<Y\>\in \BL_S$, we get $\phi\<\pi\i Y\> = \<\f\i Y\>=\<Y\>$.

For injectivity, suppose $\<\f\tilde{X}\> = \<\f\tilde{Y}\>$.  We will show $\<\pi\tilde{X}\>\leq \<\pi\tilde{Y}\>$, and injectivity follows by symmetry.  Suppose $W \in D(R)/\<\i S\>$ has $W\w \pi\tilde{Y}=0$.  Choose $\tilde{W}$ so $\pi\tilde{W}=W$.  Then $\pi(\tilde{W}\w\tilde{Y})=0$, so $\tilde{W}\w\tilde{Y}\w\i S=0$, and $0=\f(\tilde{W}\w\tilde{Y}) = \f\tilde{W}\w \f\tilde{Y}$.

By hypothesis, this implies that $0=\f\tilde{W}\w\f\tilde{X} = \f(\tilde{W}\w\tilde{X})$, so $\tilde{W}\w\tilde{X}\w\i S=0$.  This says that $0=\pi\tilde{W}\w\pi\tilde{X}=W\w\pi\tilde{X}$, and we conclude that $\phi$ is injective.

The inverse of $\phi$ is clearly given by $\phi^{-1}\<Y\> = \<\pi\i Y\>$, and both $\phi$ and $\phi^{-1}$ commute with arbitrary joins.  \end{proof}


\section{Non-noetherian rings} \label{section-nnr}

Here we will investigate maps between several graded non-Noetherian rings.  \textbf{All rings and modules in this section are graded, and objects in derived categories are bigraded.}

\begin{definition}  Fix a prime $p$.
\begin{enumerate}
   \item  For $i\geq 1$, fix integers $n_i>1$  and set
\[ \G = \frac{\zp[x_1,x_2,...]}{(x_1^{n_1}, x_2^{n_2},...)},   \;\;\L = \frac{\Fp[x_1,x_2,...]}{(x_1^{n_1}, x_2^{n_2},...)} \;\;\mbox{ and } \;\; \Pi = \frac{\bq[x_1,x_2,...]}{(x_1^{n_1}, x_2^{n_2},...)}. \]  Grade the $x_i$ so that $\G$, $\L$, and $\Pi$ are graded-connected and finitely-generated in each module degree, for example by setting $\deg(x_i) = 2^i$.  
   \item Fix $g\colon \G\to \G/p\G = \L$ to be the projection map, and fix $h\colon \G \hookrightarrow \Pi$ to be inclusion.
   \item Let $\g\colon D(\G)\to D(\L)$ and $\h\colon D(\G)\to D(\Pi)$ be the induced functors on unbounded derived categories of chain complexes of graded modules, as in Section \ref{section-BLbackground}.  Let $\ig$ and $\ih$ denote their corresponding right adjoints.
 \end{enumerate}
 \end{definition}

\begin{remark} \label{ch-cpx} Note that $\ig\L$ can be represented in $D(\G)$ by the chain complex
\[ \left(\cdots\to 0\to \G\stackrel{p}{\lra} \G\to 0\to\cdots\right).\]  Furthermore, $\ig\L$ fits into the exact triangle $\G\stackrel{p}\lra \G\to \ig\L$ in $D(\G)$.  \end{remark}

\begin{proposition} \label{fiX_X_gamma}  For all $X$ in $D(\L)$, we have $\g\ig X\cong X\oplus \Sigma X$.  Therefore  $\th(\g\ig X) = \th(X)$ and  $\<\g\ig X\> = \<X\>$.
\end{proposition} 

\begin{proof}   Using $\g \G = \L$, and Remarks \ref{Gz} and \ref{ch-cpx}, we see that $\g\ig\L$ in $D(\L)$ is 
\[ \left(\cdots\to 0\to \L\stackrel{0}{\lra} \L\to 0\to\cdots\right), \] which is just $\L\oplus \Sigma \L$.  

Next consider the small complexes, which are those in $\th(\L)$.  Let $X$ be a cofibrant representative for an arbitrary element of $\th(\L)$; then $X$ is a bounded complex of finitely-generated projective $\L$-modules.  

Since $\L$ is a local ring, projectives are free.  Thus $X$ has the form 
\[ \cdots\to \coprod_{I_3} \L\stackrel{d_2}{\lra} \coprod_{I_2} \L\stackrel{d_1}{\lra}\coprod_{I_1} \L\stackrel{d_0}{\lra}\coprod_{I_0} \L\to 0.\]

Each differential is a direct sum of maps $\L\to\L$, which we can think of as elements of $\L$.

Since every object is fibrant, $\ig X = g^*X$, and this is the complex 
\[
\UseComputerModernTips
\xymatrix{
 \cdots\ar[r]&  \coprod_{I_2}\G\ar[dr]^{\oplus p}\ar@{}[d] |{\bigoplus} \ar[r]^{\overline{d_1}} & \coprod_{I_1}\G \ar@{}[d] |{\bigoplus} \ar[dr]^{\oplus p} \ar[r]^{\overline{d_0}} & \coprod_{I_0}\G \ar@{}[d] |{\bigoplus} \ar[dr]^{\oplus p} \ar[r] & 0 \\  
 \cdots\ar[r]&  \coprod_{I_3}\G \ar[r]_{\overline{d_2}} & \coprod_{I_2}\G  \ar[r]_{\overline{d_1}}& \coprod_{I_1}\G  \ar[r]_{\overline{d_0}} & \coprod_{I_0} \G  \ar[r] & 0  .
}
\]

Here $\overline{d_i}$ is a direct sum of maps $\G\to\G$ that correspond to preimages via $g\colon\G\to\L$ of the elements of $\L$ comprising $d_i$, chosen in a compatible way.  We claim that this complex is cofibrant.  First note that it is the cofiber in $Ch(\G)$ of $W\stackrel{\oplus p}{\lra}{W}$, where $W$ is the complex
\[ \cdots\to \coprod_{I_3} \G\stackrel{\overline{d_2}}{\lra} \coprod_{I_2} \G\stackrel{\overline{d_1}}{\lra}\coprod_{I_1} \G\stackrel{\overline{d_0}}{\lra}\coprod_{I_0} \G\to 0. \]

Since $X$ is cofibrant, so is $W$.  For $X$ is an increasing union of complexes, such that the associated quotients are complexes of free $\L$-modules with zero differentials.  By replacing each $\L$ with $\G$ and each map, thought of as an element of $\L$, with its preimage via $g$ (in a compatible way), we construct $W$ as such an increasing union.  Since $W$ is cofibrant, so is $g^*X$.  

Therefore we can compute $\g(\ig X) = \g(g^*X) = g_*(g^*X)$.  But $g_*(W) = X$, and $g(p)=0$, so this gives a map $\g(\ig X) \to X\oplus \Sigma X$ that is an isomorphism.  Note that this map is functorial in $X$.

The case of a general object in $D(\L)$ follows immediately, since every object is a homotopy colimit of objects in $\th(\L)$, and $\g$ and $\ig$ are both exact and commute with coproducts. \end{proof}

This proposition allows us to apply all the results of Section \ref{rmbl-fiX-X} to the case $\g\colon D(\G)\to D(\L)$.  

Next,  we point out an important difference between $D(\G)$ and $D(\L)$.  Recall from Example \ref{example-lambda} that $\BA(\L)$ is trivial, and the module $\<I\L\>$ is a minimum nonzero Bousfield class.  The latter fact plays a significant role in~\cite{[DP]}.

Recall that, given a self-map $X\stackrel{f}{\to} X$ in any derived category $D(R)$, the homotopy colimit is called the telescope $f^{-1}X$.  More explicitly, $f^{-1}X$ is  the cofiber of the map $\coprod_{i\geq 0} X_i \stackrel{1-f}{\lra}\coprod_{i\geq 0} X_i$, where $X_i=X$ for all $i$ and the map sends each summand $X_i\to X_i\coprod X_{i+1}$ by $(1-f)(x) = (x, -f(x))$.  This is a minimal weak colimit (see e.g.~\cite[Prop.~2.2.4]{[HPS]}), so for all $n$ we have
\[ H_n(f^{-1}X)\cong \varinjlim (H_n(X)\stackrel{H(f)}{\lra} H_n(X)\to\cdots). \]

\begin{proposition} The classes $\<\ig\L\>$ and $\<\ih\Pi\>$ form a complemented pair in $\BA(\G)$.  Thus there is no minimum nonzero Bousfield class in $\BL(\G)$. \label{no-min}
\end{proposition}

\begin{proof}  Consider the self-map $\G\stackrel{p}{\lra}\G$ in $D(\G)$.  The telescope $p^{-1}\G$ is quasi-isomorphic to a module concentrated in chain degree zero, with zeroth homology the $\G$-module 
\[ \varinjlim\left(\G\stackrel{p}{\lra}\G\to\cdots\right) \cong \frac{\bq[x_1,x_2,...]}{(x_1^{n_1}, x_2^{n_2},...)} = \Pi. \]

Thus we can identify $p^{-1}\G$ with $\ih\Pi$.  As noted above, the cofiber of the map $\G\stackrel{p}{\lra}\G$ is $\ig\L$.  In this situation of a telescope and cofiber, it is well-known that $\<\ig\L\>\vee\<\ih\Pi\> = \<\G\>$; see for example~\cite[Prop.~3.6.9]{[HPS]}.

To compute $\ig\L\w\ih\Pi$, we use the chain complex description of $\ig\L$ given in Remark~\ref{ch-cpx}, and find that $\ig\L\w\ih\Pi$ is represented by 
\[ \left(\cdots\to 0\to \Pi\stackrel{p}{\lra} \Pi\to 0\to\cdots\right), \] which is zero in $D(\G)$.

Therefore the classes $\<\ig\L\>$ and $\<\ih\Pi\>$ form a nontrivial complemented pair in $\BA(\G)$.  Suppose $\<Z\>$ were a minimum nonzero Bousfield class.  Then $\<Z\>\leq \<\ig\L\>$, $\<Z\>\leq \<\ih\Pi\>$, and $\ig\L\w\ih\Pi=0$ imply $\<Z\>\w \<\ig\L\vee\ih\Pi\>=0$.  This would force $Z=0$. \end{proof}

The subcategory $\th(\ig\L)$ is a thick subcategory of compact objects in $D(\G)$.  It is clearly nonzero, and the inclusion $\th(\ig\L)\subseteq \th(\G)$ is proper.  Indeed, if $\th(\ig\L) = \th(\G)$ then we would have $\<\ig\L\>=\<\G\>$, which contradicts $\ig\L\w\ih\Pi=0$.  

\begin{definition} Let $L\colon D(\G)\to D(\G)$ be finite localization away from $\th(\ig\L)$.  Let $C$ denote the corresponding colocalization; thus for each $X$ there is an exact triangle $CX\lra X\lra LX$.
\end{definition}

See~\cite[Ch.~3]{[HPS]} or~\cite{[krause-loc-survey]} for a discussion of Bousfield localization.  Recall that we say an object $X$ is \textit{$L$-acyclic} if $L(X)=0$, and \textit{$L$-local} if it is in the essential image of $L$.

\begin{definition} The inclusion $\zp\hookrightarrow \bq$ induces a morphism $\G\to\ih\Pi$ in $D(\G)$.  Let $F$ be the fiber of this map, so $F \lra \G \lra \ih\Pi$ is an exact triangle in $D(\G)$.
\end{definition}

\begin{lemma} In $D(\G)$, the object $\Sigma F$ is quasi-isomorphic to the $\G$-module $\Pi/\G$ concentrated in degree zero. \label{sigma-f}
\end{lemma}

\begin{proof} This is a straightforward calculation using the long exact sequence in homology. \end{proof}

\begin{proposition}  \label{loc_locs} The localization functor $L$ is smashing, with $L\G = \ih\Pi$ and $C\G = F$.  It has the following acyclics and locals.
\begin{equation*} 
\begin{split}
L\mbox{-acyclics} &  = \loc(\ig\L) = \<\ih\Pi\> = \<L\G\> = \<M_g\>,\\ 
L\mbox{-locals} & = \loc(\ih\Pi) = \<\ig\L\> =\<C\G\> =\<M_h\>=\<F\>.\\ 
\end{split} 
\end{equation*}
\end{proposition}

\begin{proof} All finite localizations are smashing localizations, which means $LX = L\G\w X$.  Thus the $L$-acyclics are precisely $\<L\G\>$.  Finite localization away from $\th(\ig\L)$ means also that the $L$-acyclics are $\loc(\ig\L)$.

Next we show that the $L$-acyclics are the same as $\<\ih\Pi\>$.  Suppose $X$ is $L$-acyclic.  Then $X\in\loc (\ig\L)$.  Since $\ig\L\w\ih\Pi = 0$, this implies that $X\w\ih\Pi = 0$.  Conversely, suppose that  $X\w\ih\Pi=0$.  Then
\[ \<X\w\ig\L\> \vee \<X\w\ih\Pi\> = \<X\w \G\> = \<X\>, \] so $\<X\w\ig\L\> =\<X\>$ and $\<X\>\leq \<\ig\L\>$.  Since $\ig\L\in\loc(\ig\L)$ is $L$-acyclic, we have $L\G\w\ig\L=0$, so $L\G\w X = LX = 0$, and $X$ is $L$-acyclic.

With any smashing localization, the classes $\<L\one\>$ and $\<C\one\>$ are a complemented pair, where $\one$ is the tensor unit.  Furthermore, the $L$-locals are precisely $\<C\one\>$.  Thus in the present context, since $L\colon D(\G)\to D(\G)$ is smashing and $\<L\G\> = \<\ih\Pi\>$ is complemented by $\<\ig\L\>$, we see that the $L$-locals $\<C\G\>$ are precisely $\<\ig\L\>$.

From Remark~\ref{mf-ais} we know that $a\<\ig\L\> = \<M_g\>$.  But Proposition~\ref{no-min} shows that $\<\ig\L\>$ is complemented by $\<\ih\Pi\>$, so $a\<\ig\L\>=\<\ih\Pi\>$.  Likewise, $\<\ig\L\> = \<M_h\>$.

Next we will show that $F\w \ih\Pi=0$ in $D(\G)$.  From Lemma~\ref{sigma-f}, this is true if and only if $(\Pi/\G)\w \ih\Pi=0$ in $D(\G)$.  As in the proof of Proposition~\ref{no-min}, we can identify $\ih\Pi$ with the telescope $p^{-1}\G$, so in $D(\G)$ there is an exact triangle
\[ \bigoplus \G\stackrel{1-p}{\lra} \bigoplus \G\lra \ih \Pi. \]

Applying $(\Pi/\G)\w -$ to this, we see that $(\Pi/\G)\w\ih\Pi$ is the telescope $p^{-1}(\Pi/\G)$.  This has zero homology away from degree zero, and its degree zero homology also vanishes because $p\in \G$ so the direct limit has all zero maps.  Therefore $F\w \ih \Pi=0$, and $F$ is $L$-acyclic.  

From the triangle $F\to \G\to \ih\Pi$, we get a triangle
\[ L\G\w F \lra L\G\lra L\G\w \ih\Pi. \]

Since $L\G\w F=0$, $L\G\cong L\G\w\ih\Pi$.  Since $\ih\Pi\in \<\ig\L\>$, it is $L$-local and $L\G\w\ih\Pi\cong \ih\Pi$.  Therefore $L\G\cong \ih\Pi$.  The exact triangle $C\G\to \G\to L\G$ then forces $C\G\cong F$.

It only remains to show that the $L$-locals are given by $\loc(\ih\Pi)$.  But with a smashing localization, the $L$-locals are always form a localizing subcategory, and in addition when $T=\loc(\one)$ we always have $L$-locals = $\loc(L\one)$.  \end{proof}

The last two propositions show that $\<\ig\L\>\vee \<M_g\>=\<\G\>$ in $\BL(\G)$, so Lemma~\ref{i_BA} and Propositions~\ref{G_i_sl} and \ref{rs-lat-isomo} apply in full to this setting.

Recall that $J_g= \{\<X\>\in \BL(\G)\;|\; \g\<X\> = \0\}$.  The next theorem shows that the lattice map in Proposition~\ref{BLJ} becomes an isomorphism.

\begin{theorem} \label{blj} The functor $\g$ induces a lattice isomorphism
\[ \overline{\g}\colon \BL(\G)/J_g\to \BL(\L), \] with inverse $\ig$.
\end{theorem}

\begin{proof} Proposition~\ref{BLJ} showed that $J_g = \<M_g\>\da=\<\ih\Pi\>\da$ is a  principal ideal, and $\overline{\g}$ is a join-morphism.  We know $\ig\colon \<Y\>\mapsto [\<\ig Y\>]$ is a  join-morphism, and must only check these are inverses.

Proposition~\ref{fiX_X_gamma} shows that $\overline{\g}\ig\<X\> = \<\g\ig X\> = \<X\>$ for all $\<X\>$.

As noted earlier, we always have $\<\ig\g Y\>\leq \<Y\>$ for all $\<Y\>$.  To prove $[\ig\overline{\g}\<Y\>] =[\<Y\>]$, it remains to show that $\<\ig\g Y\>\vee\<M_g\>\geq \<Y\>\vee\<M_g\>$ for all $\<Y\>$.

So take $W$ with $W\w\ig\g Y=0$, and $W\in\<M_g\>=\<\ih\Pi\>$.  From the last proposition we get that $W$ is $L$-acyclic, so $W\w L\G=0$.  But then $W\w Y\w L\G=0$, so $W\w Y$ is also $L$-acyclic.

Now $W\w\ig\g Y=0$ implies $\g W \w \g Y = \g(W\w Y)=0$, so $W\w Y\w\ig\L=0$.  Using the last proposition, this says that $W\w Y$ is $L$-local.

Any object that is both acyclic and local with respect to a localization functor must be zero, because there are no nonzero morphisms from an acyclic to a local object.  So we conclude that $W\w Y=0$, and therefore $W\in \<Y\>\vee\<M_g\>$.  \end{proof}

Our next goal is to show that this is actually a splitting of lattices.  Towards this end, we prove some slightly more general statements.  Assume that $\cC=\loc(\one)$.  Let $l\colon\cC\to \cC$ be a smashing localization.  Define $c\one$ by the exact triangle $c\one\to \one\to l\one$.  Then for every $X\in \cC$ we have $lX\cong l\one\w X$, and an exact triangle $c\one\w X\to X \to l\one\w X$.  It follows that the $l$-acyclics are precisely $\<l\one\> = \loc(c\one)$, and the $l$-locals are $\<c\one\>=\loc(l\one)$.

Let $i_c\colon\loc(c\one)\to \cC$ and $i_l\colon\loc(l\one)\to \cC$ denote the inclusions.  

\begin{definition} Given a localizing subcategory $\cB\subseteq \cC$, let 
\[ \BL(\cC)|_\cB = \{ \<X\>\in \BL_\cC\;|\; X\in \cB\}\subseteq \BL(\cC). \]  
But note that $\<X\>=\<Y\>$ and $X\in \cB$ does not  imply  $Y\in\cB$ in general.
\end{definition}

\begin{lemma} \label{incl-gh}The inclusions $i_c$ and $i_l$ induce  join-morphisms on Bousfield lattices that preserve arbitrary joins.
\[ i_c\colon \BL(\loc(c\one))\to \BL(\cC), \mbox{ where } i_c\<X\> = \<i_c X\> = \<X\>, \mbox{ and } \]
\[ i_l\colon \BL(\loc(l\one))\to \BL(\cC), \mbox{ where } i_l\<Y\> = \<i_l Y\> = \<Y\>.\;\;\;\;\;\;\;\;\; \]
\end{lemma}

\begin{proof}  Take $X,Y\in \loc(c\one)$ such that $\<X\>\leq \<Y\>$ in $\BL(\loc(c\one))$.  Now suppose $W\in \cC$ has $W\w Y=0$.  Then $W\w c\one\w Y=0$.  But $W\w c\one \in \loc (\one\w c\one) = \loc(c\one)$, so by hypothesis we have $W\w c\one \w X=0$.  Thus $W\w X\in\<c\one\>$ is $l$-local.

Since $X\in \loc(c\one)$, it is $l$-acyclic, and $W\w X$ is also $l$-acyclic.  Any object that is both acyclic and local must be zero, so $W\w X=0$ as desired.  Therefore $i_c$ induces an order-preserving and well-defined map on Bousfield lattices.  Coproducts in both $\loc(c\one)$ and $\cC$ are given by degreewise direct sums of modules.  So $i_c$ preserves arbitrary coproducts on the object level, and thus arbitrary joins on the level of Bousfield classes.

A similar argument shows the same for $i_l$.  \end{proof}

This lemma does not generalize to arbitrary localizing subcategory inclusions, but we do have the following lemma, which is easy to prove.

\begin{lemma} \label{incl-sub} If $i\colon\cB\to \cC$ is the inclusion of a localizing tensor ideal, and induces a join-morphism on Bousfield lattices, then $\BL(\cB)\cong \BL(\cC)|_\cB$.
\end{lemma}

\begin{proposition} \label{congs} The following hold.
\[ \BL(\loc(c\one))\cong \BL(\cC)|_{\loc(c\one)} = \<c\one\>\da , \mbox{ and } \]
\[ \BL(\loc(l\one))\cong \BL(\cC)|_{\loc(l\one)} = \<l\one\>\da.\;\;\;\;\;\;\;\;\;\;\; \]
\end{proposition}

\begin{proof}  The isomorphisms on the left come from  Lemmas~\ref{incl-gh} and \ref{incl-sub}.  For the  equalities on the right, we need to know that for all $X\in \cC$ we have $\<X\>\leq \<c\one\>$ if and only if $X\in \loc(c\one)$, and $\<X\>\leq \<l\one\>$ if and only if $X\in \loc(l\one)$.  It is always the case that $X\in \loc(Y)$ implies $\<X\>\leq \<Y\>$.  If $\<X\>\leq \<l\one\>$, then since $c\one\w l\one=0$ we have $X\in \<c\one\> = \loc(l\one)$.  Similarly, if $\<X\>\leq \<c\one\>$, then $X\in \<l\one\> = \loc(c\one)$.   \end{proof}

\begin{theorem} \label{split-1} There is a lattice isomorphism
\[ \Phi\colon \BL(\cC)\stackrel{\sim}{\lra} \BL(\cC)|_{\loc(c\one)}\times  \BL(\cC)|_{\loc(l\one)}, \mbox{ where }\]
\[ \Phi\<X\>= \left(\<X\w c\one\>, \<X\w l\one\>\right).\]
The inverse is given by $\Phi'\colon (\<X\>, \<Y\>)\mapsto \<X\>\vee\<Y\>$.
\end{theorem}

\begin{proof}  Note that $X\w Z\in \loc(\one\w Z) = \loc(Z)$ for any $Z$.  It's clear that both $\Phi$ and $\Phi'$ are lattice join-morphisms.  We compute $\Phi'\Phi\<X\>$ as
\[ \<X\w c\one\>\vee\<X\w l\one\> = \<X\>\w\left(\<c\one\>\vee\<l\one\>\right) = \<X\>\w\<\one\> = \<X\>. \]

On the other hand, for $X\in\loc(c\one)=\<l\one\>$ and $Y\in\loc(l\one)=\<c\one\>$, we can compute $\Phi\Phi' (\<X\>,\<Y\>)$ as
\[ \left(\<X\vee Y\>\w\<c\one\>, \<X\vee Y\>\w\<l\one\>\right) = \left(\<X\>,\<Y\>\right), \] because $X$ is $l$-acyclic and $Y$ is $l$-local. \end{proof}

\begin{remark} Most of Lemma~\ref{incl-gh}, Proposition~\ref{congs} and Theorem~\ref{split-1} are contained, in less detail, in Proposition 6.12 and Lemma 6.13 of~\cite{[IK]}.  
\end{remark}

We now apply these general results  to our specific context.  Recall that in Proposition~\ref{loc_locs} we constructed a smashing localization functor on $D(\G).$  Let $i_g\colon \loc(\ig\L)\to D(\G)$ and $i_h\colon \loc(\ih\Pi)\to D(\G)$ denote the inclusions.  

\begin{corollary} \label{same} The functor $\g$ induces a lattice isomorphism
\[ \g\colon \BL(\loc(\ig\L))\stackrel{\sim}{\lra} \BL(\L), \mbox{ with } \<X\>\mapsto \<\g X\>. \]
\end{corollary}

\begin{proof}   This follows by combining Propositions~\ref{loc_locs} and~\ref{congs}, and Theorems~\ref{blj} and ~\ref{split-1}, and the fact that $J_g = \<\ih\Pi\>\da \cong \BL(\loc(\ih\Pi))$.

\[ \BL(\loc(\ig\L)) \cong \frac{\BL(\loc(\ig\L))\times \BL(\loc(\ih\Pi))}{0\times \BL(\loc(\ih\Pi))}\cong \BL(\G)/J_g \stackrel{\sim}{\to} \BL(\L). \]
\end{proof}

\begin{corollary} \label{split-2}The functors $\g$ and $\h$ induce a lattice isomorphism
\[ \BL({\G}) \cong \BL({\L})\times \BL({\loc(\ih\Pi)}), \mbox{ where } \]
\[ \<X\> \mapsto \left(\g\<X\>, \<X\w\ih\Pi\>\right). \]  The inverse is given by
\[ \left(\<Y\>, \<Z\>\right) \mapsto \ig\<Y\>\vee \<i_h Z\>. \]
\end{corollary}

\begin{proof}  This is an application of Theorem~\ref{split-1}, along with the observation that
\[ \g\<X\w\ig\L\> = \<\g X\>\w \<\g\ig\L\> = \<\g X\>\w \<\L\> = \<\g X\>. \]
 \end{proof}

\begin{corollary} \label{split-3} The isomorphism in Corollary \ref{split-2} induces a splitting of the distributive lattices and Boolean algebras
\[ \DL({\G}) \cong \DL({\L})\times \DL({\loc(\ih\Pi)}), \]
\[ \BA({\G}) \cong \BA({\L})\times \BA({\loc(\ih\Pi)}). \]
\end{corollary}

\begin{proof} Much of this follows from Proposition~\ref{G_i_sl}.  First consider the distributive lattice. If $\<X\>\in \DL(\G)$, then 
\[ \<X\w \ih\Pi\> \w \<X\w \ih\Pi\> = \<X\w X \w \ih\Pi\w \ih\Pi\>= \<X\w \ih\Pi\>, \] because $\ih$ maps $\DL(\Pi)$ into $\DL(\G)$.

If we consider $\<Y\>\in \DL(\L)$ and $\<Z\>\in \DL(\loc(\ih\Pi))$, then $\<\ig Y\>$ and $\<i_h Z\>$ are both in $\DL(\G)$, so their join is as well.  

Now consider the Boolean algebra.  Recall that the maximum Bousfield class of $\BL({\loc(\ih\Pi)})$ is $\<\ih\Pi\>$, and it is to this that we require complements to join.  Taking $\<X\>\in \BA(\G)$, with complement $\<X^c\>$, we compute
\[ \<X\w \ih\Pi \> \w \<X^c \w \ih\Pi\> = \0, \mbox{ and } \]
\[ \<X\w \ih\Pi\> \vee \<X^c\w \ih\Pi\> = \<\G\> \w \<\ih\Pi\> = \<\ih\Pi\>. \]

If we take  $\<Y\>\in \BA(\L)$ then Proposition~\ref{G_i_sl} implies $\<\ig Y\>\in \BA(\G)$.  For $\<Z\>\in \BA(\loc(\ih\Pi))$, with complement $\<Z^c\>$ there, one can show that $\<i_h Z\>$ is complemented in $\BL(\G)$, with complement $\<i_h (Z^c)\> \vee \<\ig \L\>$.  Therefore $\<\ig Y\>\vee\<i_h Z\>\in \BA(\G)$. \end{proof}

\begin{corollary} The Bousfield lattice of $D(\G)$ has cardinality $2^{2^{\aleph_0}}$. \label{size}
\end{corollary}

\begin{proof} Corollary B in~\cite{[DP]} states that the Bousfield lattice of $D(\L)$ has cardinality  $2^{2^{\aleph_0}}$, so $\BL({\G})$ is at least as large.  However, $\G$ is countable, so~\cite[Thm.~1.2]{[DP01]} implies that $\BL({\G})$ has cardinality at most  $2^{2^{\aleph_0}}$.
\end{proof}


\end{document}